\documentclass[10pt, a4paper, oneside]{amsart}
\usepackage{amstext, amscd, amsthm, bm, amssymb, dsfont, mathtools}
\usepackage{amsmath}

\usepackage[utf8]{inputenc}
\usepackage[T1]{fontenc}
\usepackage{lmodern}

\usepackage{cite}
\usepackage{color}
\usepackage{datetime}
\usepackage{microtype}
\usepackage[english]{babel}
\usepackage{hyperref}
\usepackage{thmtools, thm-restate}
\usepackage{tikz}
\usetikzlibrary{babel}
\usepackage{tikz-cd}
\usepackage{booktabs}
\usepackage[backgroundcolor=white, linecolor=lightgray,
bordercolor=lightgray, textsize=tiny]{todonotes}

\usepackage{cleveref}

\newcommand{\opp}{\mathrm{opp}}
\newcommand{\R}             {\mathbb R} 
    
\newcommand{\N}             {\mathbb N}
\newcommand{\Z}             {\mathbb Z}
\newcommand{\Q}             {\mathbb Q}

\newcommand{\cat}[1]{\pmb{\mathrm{#1}}}
\newcommand{\Rder}{\cat{R}}
\newcommand{\Lder}{\cat{L}}
\newcommand{\Dder}{\cat{D}}

\usetikzlibrary{arrows,calc,matrix}

\usepackage{tikzpfeile}

\makeatletter
\newcommand{\osetwv}[2]{
  {\mathop{#2}\limits^{\vbox to -4\ex@{\kern-\tw@\ex@
   \hbox{\hspace{-0.2em}\scriptsize #1}\vss}}}}
\makeatother

\makeatletter
\newcommand{\oset}[2]{
  {\mathop{#2}\limits^{\vbox to -5\ex@{\kern-\tw@\ex@
   \hbox{\hspace{-0.25em}\scriptsize$#1$}\vss}}}}
\makeatother

\crefname{sequence}{sequence}{sequence}
\let\realcite\cite
\renewcommand{\ref}[1]{\cref{#1}}
\renewcommand{\cite}[2][]{\hspace{-0.001ex}{\realcite[#1]{#2}}}
\newcommand{\citenospec}[1]{\hspace{-0.001ex}{\realcite{#1}}}

\DeclarePairedDelimiter\abs{\lvert}{\rvert}%
\DeclarePairedDelimiter\norm{\lVert}{\rVert}%

\makeatletter
\let\oldabs\abs
\def\abs{\@ifstar{\oldabs}{\oldabs*}}
\let\oldnorm\norm
\def\norm{\@ifstar{\oldnorm}{\oldnorm*}}
\makeatother

\DeclareMathOperator{\Hom}{Hom}

\DeclareMathOperator{\Ext}{Ext}
\DeclareMathOperator{\Tor}{Tor}
\DeclareMathOperator{\tor}{tor}

\DeclareMathOperator{\id}{id}
\DeclareMathOperator{\RHom}{\cat{R}Hom}
\DeclareMathOperator{\E}{E}

\newcommand{\RMod}[1]{#1\textrm{-}\cat{Mod}}

\declaretheorem[numberwithin=section]{theorem}

\newenvironment{thm*}[1]{\par{\bfseries #1}.\itshape}{}
\newenvironment{thmcite*}[2]{\par{\bfseries #1} (#2).\itshape}{}

\declaretheorem[sibling=theorem]{corollary}
\declaretheorem[sibling=theorem]{lemma}

\declaretheorem[sibling=theorem]{proposition}

\declaretheorem[style=remark, sibling=theorem]{remark}

\declaretheorem[style=remark, sibling=theorem]{example}

\declaretheorem[style=definition, sibling=theorem]{definition}

\newcommand{\rom}[1]{\uppercase\expandafter{\romannumeral #1\relax}}

\address{Mathematisches Institut\\Universität Heidelberg\\Im
  Neuenheimer Feld 205\\D-69120 Heidelberg}\author{Oliver Thomas}

\email[O.\,Thomas]{othomas@mathi.uni-heidelberg.de}
\urladdr[O.\,Thomas]{https://www.mathi.uni-heidelberg.de/~othomas/}

\author{Otmar Venjakob}
\email[O.\,Venjakob]{otmar@mathi.uni-heidelberg.de}
\urladdr[O.\,Venjakob]{https://www.mathi.uni-heidelberg.de/~otmar/}

\title[Spectral Sequences for Iwasawa Adjoints]{On Spectral Sequences for Iwasawa Adjoints à la Jannsen for
  Families}

\begin{document}
\maketitle

\begin{abstract}
  In \citenospec{MR1097615} several spectral sequences for (global and
  local) Iwasawa modules over (not necessarily commutative) Iwasawa
  algebras (mainly of $p$-adic Lie groups) over $\Z_p$ are
  established, which are very useful for determining certain
  properties of such modules in arithmetic applications. Slight
  generalizations of said results can be found in
  \citenospec{MR2333680} (for abelian groups and more general
  coefficient rings), \citenospec{MR1924402} (for products of not
  necessarily abelian groups, but with $\Z_p$-coefficients), and
  \citenospec{MR3084561}. Unfortunately, some of Jannsen's spectral
  sequences for families of representations as coefficients for
  (local) Iwasawa cohomology are still missing. We explain and follow
  the philosophy that all these spectral sequences are consequences or
  analogues of local cohomology and duality à la Grothendieck (and
  Tate for duality groups).
\end{abstract}

\tableofcontents

\section{Introduction}
\label{sec:introduction}

Let $\mathcal O$ be a complete discrete valuation ring with
uniformising element $\pi$ and finite residue field. Consider
furthermore an $\mathcal O$-algebra $R$, which is a complete
Noetherian local ring with maximal ideal $\mathfrak m$, of dimension
$d$ and finite residue field. We are mainly interested in the case of
a ring of formal power series $R=\mathcal O[[X_1,\dots, X_t]]$ in $t$
variables, which is a complete regular local ring of dimension
$d=t+1$. We now have a number of dualities at hand.

First, there is Matlis duality: Denote with $\mathcal E$ an injective
hull of $R/\mathfrak m$ as an $R$-module. Then
$T=\Hom_R(-,\mathcal E)$ induces a contravariant involutive
equivalence between Noetherian and Artinian $R$-modules akin to
Pontryagin duality.

Second, there is local duality: If
$\Rder\Gamma_{\underline{\mathfrak m}}$ denotes the right derivation
of $ M \mapsto \varinjlim_k \Hom_R(R/\mathfrak m^k, M)$ in the derived
category of $R$-modules, then
\[ \Rder\Gamma_{\underline{\mathfrak m}} \cong [-d] \circ T \circ
  \RHom_R(-, R)\] on coherent $R$-modules.

Third, there is Koszul duality: The complex
$\Rder\Gamma_{\underline{\mathfrak m}}$ can be computed by means of
Koszul complexes $K^\bullet$ which are self-dual: $K^\bullet =
\Hom_R(K^\bullet, R)[d]$.

Finally, there is Tate duality: Let $G$ be a pro-$p$ duality group of
dimension $s$. Then for discrete $G$-modules $A$ we have
$H^i(G, \Hom(A,I)) \cong H^{s-i}(G,A)^*$ for a dualizing module $I$.

Consider $\Lambda_R(G)=\varprojlim_U R[G/U]$ where $U$ runs through
the open normal subgroups of $G$. It is well known that
$\Lambda_R(\Z_p^s)\cong R[[Y_1,\dots,Y_s]]$ and indeed
$R=\Lambda_{\mathcal O}(\Z_p^r)$. The maximal ideal of $\Lambda_R(G)$
is now generated by the regular sequence
$(\pi, X_1, \dots, X_t, Y_1, \dots, Y_s)$ and no matter how we split
up this regular sequence into two, they will remain regular. The
Koszul complex then gives rise to a number of interesting spectral
sequences and these should (at least morally) recover the spectral
sequences
\begin{equation} \Tor_n^{\Z_p}(D_m(M^\vee), \Q_p/\Z_p) \Longrightarrow
  \Ext_{\Lambda_{\Z_p}(G)}^{n+m}(M, \Lambda_{\Z_p}(G))^\vee
  \label{eqn:jannsen-pont-of-iw-adjoints-is-tor-tate-pont}
\end{equation}
and
\begin{equation} \varinjlim_k D_n(\Tor^{\Z_p}_m( \Z_p/p^k, M)^\vee \Longrightarrow
  \Ext_{\Lambda_{\Z_p}(G)}^{n+m}(M, \Lambda_{\Z_p}(G))^\vee),
  \label{eqn:jannsen-pont-of-iw-adjoints-is-colim-of-tate-of-pont-of-loc-coh}
\end{equation} which
show up in Jannsen's proof of \cite[2.1 and 2.2]{MR1097615}. The
functors $D_n$ stem from Tate's spectral sequence and are a corner
stone in the theory of duality groups.

We will show in \cref{sec:iwasawa-adjoints} that these spectral sequences (and
many more) are consequences of the four duality principles laid out above. This
also allows us to generalize Jannsen's spectral sequences to more general
coefficients. For example, generalisations of
\cref{eqn:jannsen-pont-of-iw-adjoints-is-tor-tate-pont} and
\cref{eqn:jannsen-pont-of-iw-adjoints-is-colim-of-tate-of-pont-of-loc-coh} are
subject of \cref{prop:pont-of-iw-adjoints-is-tor-tate-pont} and of
\cref{prop:pont-of-iw-adjoints-is-colim-of-tate-of-pont-of-loc-coh}
respectively. While another spectral sequence for Iwasawa adjoints has already
been generalized to more general coefficients
(cf.~\cref{thm:lim-sharifi-spec-seq-local}), the generalizations of the
aforementioned spectral sequences are missing in the literature. We can even
generalize an explicit calculation of Iwasawa adjoints (cf.~\cite[corollary
2.6]{MR1097615}, \cite[(5.4.14)]{MR2392026}) in
\cref{prop:dual-of-iw-adj-is-twist-of-local-coh-of-coeff}.

Furthermore, we generalize Venjakob's result on local duality for
Iwasawa algebras (\cite[theorem 5.6]{MR1924402}) to more general
coefficients (cf.~\cref{thm:local-duality-for-iw-alg}).  As an
application we determine the torsion submodule of local Iwasawa
cohomology generalizing a result of Perrin-Riou in the case $R=\Z_p$
in \cref{thm:torsion-in-local-iw-coh-for-poinc-grps}.

\section{Conventions}

A \emph{ring} will always be unitary and associative, but not
necessarily commutative. If not explicitly stated otherwise, “module”
means left-module, “Noetherien” means left-Noetherien etc.

We will furthermore use the language of derived categories. If $\cat{A}$ is an
abelian category, we denote with $\Dder(\cat{A})$ the derived category of
unbounded complexes, with $\Dder^+(\cat{A})$ the derived category of complexes
bounded below, with $\Dder^-(\cat{A})$ the derived category of complexes bounded
above and with $\Dder^b(\cat{A})$ the derived category of bounded complexes.

As we simultaneously have to deal with left- and right-exact functors, both
covariant and contravariant, recovering spectral sequences from isomorphisms in
the derived category is a bit of a hassle regarding the indices. Suppose that
$\cat{A}$ has enough injectives and projectives and that $M$ is a (suitably
bounded) complex of objects of $A$. Then for a covariant functor $F\colon
\cat{A}\ra\cat{A}$ we set $\Rder F=F(Q)$ and $\Lder F=F(P)$ with $Q$ a complex
of injective objects, quasi-isomorphic to $A$ and $P$ a complex of projectives,
quasi-isomorphic to $A$. If $F$ is contravariant, we set $\Lder F(A)=F(Q)$ and
$\Rder F(A)=F(P)$. For indices, this implies the following: Assume that $A$ is
concentrated in degree zero. Then for $F$ covariant, $\Rder F(A)$ has
non-vanishing cohomology at most in non-negative degrees and $\Lder F(A)$ at
most in non-positive degrees. For $F$ contravariant, it's exactly the other way
around. We set $\Lder^qF(A)=H^q(\Lder F(A))$ and $\Rder^qF (A)=H^q(\Rder F(A))$.
Note that with these conventions
\begin{itemize}
\item $\Rder^p (-)^G (A) = H^p(G,A)$
\item $\Lder^q (-\otimes M)(A) = \Tor_{-q}(A,M)$
\item $\Rder^p \Hom(-, M)(A) = \Ext^{p}(A,M)$
\item $\Lder^q (\varinjlim_U (-^U)^*)(A) = \varinjlim_U H^{-q}(U,A)^*$
\end{itemize}

If $F\colon \cat{A}\ra\cat{A}$ is exact, then $F$ maps quasi-isomorphic
complexes to quasi-isomorphic complexes. Its derivation $\Rder F$ (or $\Lder F$)
is then given by simply applying $F$ and we will make no distinction between
$F\colon \cat{A}\ra\cat{A}$ and $\Rder F\colon \Dder(\cat{A})\ra \Dder(\cat{A})$
in this case.

For every integer $d\in\Z$ we have a \emph{shift operator} $[d]$, so that for
complexes $C$ and $n\in\Z$ the following holds: \[ ([d](C))^n = C^{n+d}.\] We
will at times write $C[d]$ instead of $[d](C)$. Note that although we
occasionally cite \citenospec{MR1269324}, we deviate from Weibel's conventions
in this regard: Our $[d]$ is Weibel's $[-d]$. We furthermore set
$\Hom(C^\bullet, D^\bullet)$ to be the complex with entries $\Hom(C^\bullet,
D^\bullet)^i = \bigoplus_{n\in\Z} \Hom(C^k, D^{k+1})$. Sign conventions won't
matter in this paper.

Recall that an $R$-module is called \emph{coherent} if it is finitely
generated and all of its finitely generated submodules are finitely
presented. The natural functor from coherent $R$-modules to
$R$-modules then induces equivalences
$\Dder^*(\cat{Coh}(R))\cong \Dder^*_c(\RMod{R})$ for $*\in\{+,b\}$
where subscript “c” means complexes with coherent cohomology. If $R$
is Noetherian, then the notions of coherent, finitely generated and
Noetherian modules coincide.

\section{A few facts on $R$-modules}

\subsection{Noncommutative rings}

Let $R$ be a ring. The intersection of all maximal left ideals
coincides with the intersection of all maximal right ideals and is
called the \emph{Jacobson radical} of $R$ and is hence a two-sided
ideal, denoted by $J(R)$.  For every $r\in J(R)$ the element $1-r$
then has both a left and a right inverse and the following form of
Nakayama's lemma holds, cf.~e.\,g.~\cite[4.22]{MR1125071}.

\begin{lemma}[Azumaya-Krull-Nakayama]
  Let $M\in\RMod{R}$ be a finitely generated $R$-module. If
  $J(R)M=M,$ then $M=0$.
\end{lemma}

Recall that a ring is called local if it has a unique maximal
left-ideal. This unique left ideal is then also the ring's unique
maximal right-ideal and the group of two-sided units is the complement
of this maximal ideal.

The following is well known and gives rise to the notion of “finitely
presented” (or “compact”) objects in arbitrary categories.

\begin{lemma}
  \label{lm:hom2-commutes-with-dir-lim-iff-m-fin-pres}
  Let $R$ be a ring and $M$ an $R$-module. Then $\Hom_R(M, -)$
  commutes with all direct limits if and only if $M$ is finitely
  presented.
\end{lemma}

If $R$ is Noetherien, this isomorphism extends to higher Ext-groups.

\begin{proposition}\label{lm:ext-commutes-with-dir-lim-if-r-noeth-m-fg}
  Let $R$ be a Noetherian ring, $M$ a finitely generated
  $R$-module, and $(N_i)_i$ a direct system of $R$-modules. Then
  \[ \Ext^q_R(M, \varinjlim_i N_i) \cong \varinjlim_i \Ext^q_R(M,
    N_i).\]
  \begin{proof}
    As $R$ is Noetherian, there exists a resolution of $M$ of finitely
    generated projective $R$-Modules. As $\varinjlim$ commutes with
    homology,
    \cref{lm:hom2-commutes-with-dir-lim-iff-m-fin-pres} yields
    the result.
  \end{proof}
\end{proposition}

\begin{remark}
  Recall the following subtleties: Let $R,S,T$ be rings, $N$ a
  $S$-$R$-bimodule and $P$ a $S$-$T$-bimodule. Then $\Hom_S(N,P)$ has
  the natural structure of an $R$-$T$-bimodule via $(rf)(n) = f(nr)$
  and $(ft)(n)=f(n)t$.

  Furthermore let $M$ be a $R$-left-module. Then canonically
  \[ \Hom_R(M, \Hom_S(N,P)) = \Hom_S(N\otimes_R M, P)\] as $T$-modules.
\end{remark}

\begin{lemma}\label{lm:hom-of-flat-inj-is-inj}
  If $P$ is an $R$-$R$-bimodule that is flat as an $R$-right-module
  and $Q$ an injective $R$-module, then $\Hom_R(P,Q)$ is again an
  injective $R$-left-module.
  \begin{proof}
    $\Hom_{R}(-,\Hom_R(P,Q)) = \Hom_R(-,Q) \circ (P\otimes_R -)$ is a
    composition of exact functors and hence exact.
  \end{proof}
\end{lemma}

\begin{lemma}\label{lm:tor-m-n-zero-then-ext-m-hom-n-q-zero}
  Let $N$ be an $R$-$R$-bimodule and $M$ an $R$-modules. Then
  \[\Tor_q^R(N,M)=0 \text{ for } 1\leq q \leq n\] if and only if
  \[\Ext^q_R(M,\Hom_R(N,Q))=0 \text{ for } 1 \leq q \leq n\] and all injective
  $R$-left-modules $Q$.
  \begin{proof}  
    Let \[0 \ra K \ra P \ra N \ra 0\] be an exact sequence of
    $R$-$R$-bimodules with $P$ flat from the right. First of all,
    applying $h_Q=\Hom_R(-,Q)$ yields the exact sequence
    \begin{equation}0 \ra h_Q(N) \ra h_Q(P) \ra h_Q(K) \ra 0
      \label[sequence]{seq:dual-of-resolution-of-n}
    \end{equation}
    as $Q$ is injective. By \cref{lm:hom-of-flat-inj-is-inj} $h_Q(P)$ is
    injective, so $\Ext_R^q(M,h_Q(P))=0$ for $q\geq 1$ and hence
    $\Ext^{q+1}_R(-,h_Q(N)) = \Ext^q_R(-,h_Q(K))$ for $q\geq 1$. Similarly
    $\Tor_{q+1}^R(N,-)=\Tor_q^R(K,-)$ for $q \geq 1$. Hence we can by dimension
    shifting reduce to the case $n=q=1$.

    Tensoring above exact sequence with $M$ from the right yields the
    exact sequence
    \[0 \ra \Tor_1^R(N,M) \ra K\otimes M \ra P \otimes M \ra
      N\otimes M\ra 0\] and by injectivity of $Q$ the exact sequence
    \[0 \ra h_Q(N\otimes M) \ra h_Q(P\otimes M) \ra
      h_Q(K\otimes M) \ra h_Q(\Tor_1^R(N,M)) \ra 0.\]

    We can also apply $\Ext^\bullet_R(M,-)$ to the exact
    \cref{seq:dual-of-resolution-of-n}. By Tensor-Hom-adjointness and
    \cref{lm:hom-of-flat-inj-is-inj}, this implies that \[h_Q(\Tor_1^R(N,M)) =
      \Ext^1_R(M, h_Q(N)),\] so if $\Tor_1^R(N,M)=0$, then clearly $\Ext^1_R(M,
    h_Q(N))=0$ for all injective $Q$. If on the other hand $\Ext^1_R(M,
    h_Q(N))=0$ for all injective $Q$, then $h_Q(\Tor_1^R(N,M))=0$ for all
    injective $Q$ and hence also $\Tor_1^R(N,M)=0$ as the category of
    $R$-modules has sufficiently many injectives.
  \end{proof}
\end{lemma}

\begin{definition}
  Let $R$ be a ring. A sequence $(r_1,\dots,r_d)$ of central elements
  in $R$ is called \emph{regular}, if for each $i$ the residue class
  of $r_{i+1}$ in $R/(r_1,\dots,r_i)$ is not a zero-divisor.
\end{definition}

\begin{proposition}\label{prop:tor-vanishes-regular-sequences}
  Let $R$ be a ring and
  $(r_1,\dots,r_k,s_1,\dots,s_l)$ such a regular sequence that the
  sequence $(s_1,\dots,s_l)$ is itself regular. Let $I=(r_i)_i$ and
  $J=(s_i)_i$ be the ideals generated by the first and second part of
  the regular sequence. Then for all $q\geq 1$
  \[ \Tor^R_q(R/I, R/J) = 0\]
  and
    \[\varprojlim_n \Tor^R_q(R/I^n, R/J^n)=0.\]
  \begin{proof}
    Let us first show that $\Tor^R_1(R/I, R/J) = 0$ and then reduce to
    this case by induction on $l$. Consider the exact sequence of
    $R$-modules \[0 \ra J \ra R \ra R/J \ra 0\]
    and apply $\Tor^R_\bullet(R/I,-)$. As $R$ is flat,
    \[\Tor^R_1(R/I, R/J) = \ker \left(R/I\otimes_R J \ra R/I\right) = \frac{I\cap
        J}{IJ}.\]

    We argue by induction on $l$ that this is zero: For $l=1$,
    $x\in I\cap J$ implies $x=\lambda s_1 = s_1\lambda \in I$, so
    $\lambda=0$ in $R/I$, so $\lambda\in I$ and hence $x\in IJ$.
    Denote with $J'$ the ideal generated by $s_1,\dots,s_{l-1}$. By
    induction $I\cap J'=IJ'$. Let $x\in I\cap J = I\cap (J'+s_lR)$, so
    $x= a+s_lb$ with $a\in J'$ and $b\in R$. Clearly $s_l b=0$ in
    $R/(I+J')$, so $b\in I+J'$ by regularity of the sequence and hence
    \[I\cap J = I\cap (J'+s_lR) = I \cap (J' + s_l(I+J')) = I \cap (J'
      + s_lI) = I\cap J' + I\cap s_l I\] as $s_lI\subseteq I$. Now
    $I\cap J' + I\cap s_l I = IJ' +s_l I = I(J'+s_lR) = IJ$, and this
    was to be shown.

    We now argue by induction on $l$ that $\Tor^R_q(R/I,R/J)=0$ for
    all $q>0$. For $l=1$ we have the free resolution
    \[0 \ra R \oset{s_1}{\ra} R \ra R/J \ra 0\]
    hence $\Tor^R_q(R/I, R/J)=0$ for $q>1$ and for $q=1$ by what we
    saw above. Let $J'$ be again the ideal generated by
    $s_1,\dots,s_{l-1}$. By induction we can assume that all
    $\Tor_q^R(R/I,R/J')$ vanish. Consider the sequence
    \[ 0 \ra R/J' \oset{s_l}{\ra} R/J' \ra R/J \ra 0,\]
    which is exact as the subsequence $(s_1,\dots,s_l)$ is
    regular. Applying $\Tor^R_\bullet(R/I,-)$ shows that
    $\Tor^R_q(R/I,R/J)=0$ for $q>1$ by induction hypothesis — and by
    what we saw above also for $q=1$.

    It is easy to see that the sequence
    $(r_1^m, \dots, r_k^m,s_1^m,\dots,s_l^m)$ is again regular
    (cf.~e.\,g.\cite[theorem 16.1]{MR879273}). Denote the ideal
    generated by $(r_i^m)_i$ and $(s_i^m)_i$ by $I^{(m)}$ and
    $J^{(m)}$ respectively. Let $m'=\max\{k,l\}$. Clearly
    $I^{m'm}\subseteq I^{(m)}\subseteq I^m$ and the same is true for
    $J$, hence the natural map
    \[\Tor^R_q(R/I^{m'm}, R/J^{m'm}) \ra \Tor^R_q(R/I^m, R/J^m)\]
    factors through $\Tor^R_q(R/I^{(m)}, R/J^{(m)})$, which is zero.
  \end{proof}
\end{proposition}

\subsection{The Koszul Complex}

Recall the following (cf.~e.\,g.\cite[section 4.5]{MR1269324}).

\begin{definition}
  Denote by \[ K_\bullet(x) = 0 \ra R \oset{x}{\ra} R \ra 0\] the
  chain complex (i.\,e., the degree decreases to the right)
  concentrated in degrees one and zero for a ring $R$ and a central
  element $x\in Z(R)$. For central elements $x_1,\dots,x_d$ the complex
  \[K_\bullet(x_1,\dots,x_d) = K_\bullet(x_1)\otimes\dots\otimes
    K_\bullet(x_d)\] is called the \emph{Koszul complex attached to
    $x_1,\dots,x_d$}. We will also consider the cochain complex
  $K^\bullet(x_1,\dots,x_d)$ with entries
  $K^p(x_1,\dots,x_d)=K_{-p}(x_1,\dots,x_d)$.
\end{definition}

\begin{remark}
  While this definition is certainly elegant, a more down to earth
  description is given as follows: $K_{p}(x_1,\dots,x_d)$ is the free
  $R$-module generated by the symbols
  $e_{i_1}\wedge\dots\wedge e_{i_p}$ with $i_1<\dots<i_p$ with
  differential
  \[d(e_{i_1}\wedge\dots\wedge e_{i_p}) = \sum_{k=1}^p (-1)^{k+1}
    x_{i_k} e_{i_1}\wedge\dots\wedge\widehat{e_{i_k}}\wedge\dots\wedge
    e_{i_p}.\] This description emphasizes the importance of using
  central elements $(x_i)_i$.
\end{remark}

\begin{remark}
  The importance of the Koszul complex for our purposes stems from the
  following fact: If $x_1,\dots,x_d$ is a regular sequence, then
  $K_\bullet(x_1,\dots,x_d)$ is a free resolution of
  $R/(x_1,\dots,x_d)$, cf.~e.\,g.\cite[corollary 4.5.5]{MR1269324}.
\end{remark}

\begin{proposition}\label{prop:koszul-poincare-duality}
  The complex $K_\bullet=K_\bullet(x_1,\dots,x_d)$ is isomorphic to
  the complex
  \[ 0 \ra \Hom_R(K_0, R) \ra \dots \ra \Hom_R(K_{d}, R) \ra 0,\]
  where $\Hom_R(K_0,R)$ is in degree $d$ and $\Hom_R(K_d,R)$ in degree
  zero. Analogously \[K^\bullet \cong \Hom_R(K^\bullet, R)[d].\]
  \begin{proof}
    We have to describe isomorphisms $K_{p} \cong \Hom_R(K_{d-p},R)$
    such that all diagrams \[\begin{tikzcd} K_{p}
        \arrow{r}{d}\arrow{d}{\cong}& K_{p-1} \arrow{d}{\cong}\\
        \Hom_R(K_{d-p},R)\arrow{r}{d^*}&
        \Hom_R(K_{d-p+1},R)\end{tikzcd} \] commute. Consider the map
    \[ e_{i_1}\wedge\dots\wedge e_{i_p} \mapsto \left( e_{j_1} \wedge
        \dots \wedge e_{j_{d-p}} \mapsto \mathrm{sgn}(\sigma)
      \right) \] where $\mathrm{sgn}(\sigma)$ is the sign of the
    permutation
    $\sigma(1)=i_1, \dots, \sigma(p)=i_p, \sigma(p+1)=j_1,\dots,
    \sigma(d)=j_{d-p}$, i.\,e., in the exterior algebra
    $\bigwedge^dR^d$ we have
    \[e_{i_1}\wedge\dots\wedge e_{i_p} \wedge e_{j_1}\wedge \dots
      \wedge e_{j_{d-p}} = \mathrm{sgn}(\sigma) e_1\wedge\dots\wedge
      e_d.\] (Note that $\mathrm{sgn}(\sigma)=0$ if $\sigma$ is not
    bijective.) It is then easy to verify that the diagram above does
    indeed commute: Identifying $R$ with $\bigwedge^dR^d$, all but at
    most one summand vanishes in the ensuing calculation and the
    difference in sign is precisely the difference in the
    permutations.
  \end{proof}
\end{proposition}

\begin{proposition}\label{prop:duality-rhom-tor-for-reg-seq-in-der-cat}
  Let $R$ be a ring and $x_1,\dots,x_d$ a regular sequence of central
  elements in $R$. Then in the bounded derived category of $R$-modules
  \[ [d]\circ \RHom_R(R/(x_1,\dots,x_n), -) \cong R/(x_1,\dots,x_n)
    \otimes^{\cat{L}}_R -.\]
  \begin{proof}
    Denote with $K_\bullet$ the Koszul (chain) complex
    $K_\bullet(x_1,\dots,x_d)$ (concentrated in degrees $d,d-1,\dots, 0$) and
    with $K^\bullet$ the Koszul (cochain) complex (concentrated in degrees $-d,
    -d+1, \dots, 0$).

    As $x_1,\dots,x_d$ form a regular sequence, $K^\bullet$ is a free
    resolution of $R/(x_1,\dots,x_n)$ and hence allows us to calculate
    the derived functors as follows.
    \begin{align*}
      \RHom_R(R/(x_1,\dots,x_n), M)[d] &=
                                        \Hom_R(K^\bullet,M)[d]\\
                                      &\cong \Hom_R(K^\bullet[-d],R) \otimes_R M\\
                                      &\cong K^\bullet \otimes_R M\\
                                      &= R/(x_1,\dots,x_d) \otimes^{\cat{L}}_R M,
    \end{align*}
    with the crucial isomorphisms being due to the fact that
    $K^\bullet$ is a complex of free modules and
    \cref{prop:koszul-poincare-duality}. It is clear that these
    isomorphisms are functorial in $M$.
  \end{proof}
\end{proposition}

\begin{corollary}\label{prop:rhom-r-mod-regular-commutes-with-matlis}
  Let $R$ be a commutative ring, $x_1,\dots,x_d$ a regular sequence in
  $R$ and $T=\Hom_R(-,Q)$ with $Q$ injective. Then
  \[ \RHom_R(R/(x_1,\dots,x_d), -) \circ T = T \circ [d] \circ
    \RHom_R(R/(x_1,\dots,x_d), -) \] on the derived category of
  $R$-modules.
  \begin{proof}
    The functor $T$ is exact and
    \[\Hom_R(R/(x_1,\dots,x_d),-)\circ T = T \circ
      (R/(x_1,\dots,x_d)\otimes_R -)\] by adjointness. Hence
    also
    \[\RHom_R(R/(x_1,\dots,x_d),-)\circ T = T \circ
      (R/(x_1,\dots,x_d)\otimes_R^{\Lder} -).\] By
    \cref{prop:duality-rhom-tor-for-reg-seq-in-der-cat}, this is just $T \circ
    [d] \circ \RHom_R(R/(x_1,\dots,x_d), -).$
  \end{proof}
\end{corollary}

\begin{corollary}\label{prop:duality-ext-tor-for-reg-seq-classical}
  Let $R$ be a commutative ring and $x_1,\dots,x_d$ a regular sequence
  in $R$. Let further $M$ be an $R$-module. Then
  \[\Ext^{d-p}_R(R/(x_1,\dots,x_d), M)=\Tor^R_{p}(R/(x_1,\dots,x_d),
  M).\]
  \begin{proof}
    This is just \cref{prop:duality-rhom-tor-for-reg-seq-in-der-cat},
    taking extra care of the indices:
    \begin{align*}\Ext^{d-p}_R(R/(x_1,\dots,x_d), M) &=
    \cat{R}^{d-p}\Hom_R(R/(x_1,\dots,x_d), M)\\ &=
    H^{-p}\RHom_R(R/(x_1,\dots,x_d)[-d], M) \\&= H^{-p}(R/(x_1,\dots,x_d)
    \otimes^{\cat{L}}_R M) \\ &= \Tor^R_{p}(R/(x_1,\dots,x_d),M).\end{align*}
  \end{proof}
\end{corollary}

\subsection{Local Cohomology}

\begin{definition}
  Let $R$ be a ring and $\underline{J} = (J_n)_{n\in \N}$ a
  decreasing sequence of two-sided ideals. (The classical example is
  to take a two-sided ideal $J$ and set $\underline{J}=(J^n)_n$.) For
  an $R$-left-module $M$ set \[ \Gamma_{\underline{J}}(M) = \{m\in M\mid
    J_n m = 0 \text{ for some $n$}\}.\]

  It is clear that $\Gamma_{\underline{J}}$ is a left-exact functor
  with values in $R\cat{\text{-}Mod}$. Denote its right-derived
  functor in the derived category $\Dder^+(R\cat{\text{-}Mod})$ by
  $\Rder\Gamma_{\underline{J}}$.
\end{definition}

\begin{remark}\label{rem:different-expressions-of-loc-coh}
  $\Gamma_{\underline{J}} = \varinjlim_n \Hom_R(R/J_n,-)$, so
  $\Rder\Gamma_{\underline{J}} = \varinjlim_n
  \RHom_R(R/J_n,-)$ and $\Rder^q\Gamma_{\underline{J}} = \varinjlim_n
  \Ext_R^q(R/J_n,-)$.
\end{remark}


\begin{lemma}\label{prop:right-adj-preserves-inj-obj-if-left-adj-exact}
  Let $\cat{A},\cat{B}$ be abelian categories, with additive functors
  $L\colon \cat{A}\ra\cat{B}$ left adjoint to
  $R\colon \cat{B}\ra\cat{A}$. If $L$ is exact, $R$ preserves
  injective objects.
  \begin{proof}
    This is well known, cf.~e.\,g.\cite[proposition 2.3.10]{MR1269324}.
  \end{proof}
\end{lemma}

\begin{remark}\label{prop:loc-coh-indep-of-base}
  Let $\varphi\colon R\ra S$ be a homomorphism between unitary rings. Let
  $\underline{J}$ be decreasing sequence of two-sided ideals in $R$ and denote
  with $\underline{J}S$ the induced sequence of two-sided ideals in $S$. If
  $\varphi(R)$ lies in the centre of $S$, then $\Gamma_{\underline{J}} \circ
  \mathrm{res}_\varphi=\mathrm{res}_\varphi\circ \Gamma_{\underline{J}S}$. If
  furthermore injective $S$-modules are also injective as $R$-modules, e.\,g.,
  if $S$ is a flat $R$-module via
  \cref{prop:right-adj-preserves-inj-obj-if-left-adj-exact}, then
  $\Rder\Gamma_{\underline{J}} \circ
  \mathrm{res}_\varphi=\mathrm{res}_\varphi\circ\Rder\Gamma_{\underline{J}S}$.
  Local cohomology is thus independent of the base ring for flat extensions and
  we will omit $\mathrm{res}_\varphi$ and the distinction between
  $\underline{J}S$ and $\underline{J}$ in the future. Note especially that if
  $R$ is complete, then $R[[G]]=R[G]^\wedge$ is a flat $R$-module.
\end{remark}

\begin{proposition}\label{prop:loc-coh-is-max-fin-submod}
  Let $R$ be a Noetherian local ring with maximal ideal $\mathfrak m$
  and finite residue field. Let $M$ be a finitely generated
  $R$-module. Then $\Gamma_{\underline{\mathfrak m}}(M)$ is the
  maximal finite submodule of $M$.
  \begin{proof}
    Denote with $T$ the maximal finite submodule of $M$ (which exists
    as $M$ is Noetherian). By Nakayama there exists a $k\in\N$ with
    $\mathfrak m^kT=0$ and hence
    $T\subseteq\Gamma_{\underline{\mathfrak m}}(M)$. Conversely
    $R/\mathfrak m^k$ is a finite ring for each $k$, hence $Rm$ is a
    finite module for each $m\in\Gamma_{\underline{\mathfrak m}}(M)$
    and is thus contained in $T$.
  \end{proof}
\end{proposition}

\begin{proposition}\label{prop:local-coh-commutes-with-direct-limits}
  If $R$ is a Noetherian ring and $\underline{J}$ a decreasing
  sequence of ideals, then $\Rder\Gamma_{\underline{J}}$ and
  $\Rder^q\Gamma_{\underline{J}}$ commute with direct limits.
  \begin{proof}
    This is just
    \cref{lm:ext-commutes-with-dir-lim-if-r-noeth-m-fg}, as for Noetherian
    rings, direct limits of injective modules are again injective.
  \end{proof}
\end{proposition}

\begin{definition}
  For $\underline{I}$ and $\underline{J}$ decreasing sequences of
  two-sided ideals of a ring $R$ set
  $(\underline{I} + \underline{J})_n=I_n+J_n$.
\end{definition}

\begin{remark}
  If $I$ and $J$ are two-sided ideals of a ring $R$, then
  generally $\underline{I}+\underline{J}\neq \underline{I+J}$. But as
  these two families are cofinal,
  $\Gamma_{\underline{I}+\underline{J}} = \Gamma_{\underline{I+J}}$.
\end{remark}

\begin{remark}
  Clearly
  $\Gamma_{\underline{I}+\underline{J}}=\Gamma_{\underline{I}}\circ\Gamma_{\underline{J}}$,
  but regrettably
  \[\Rder\Gamma_{\underline{I}+\underline{J}} = \Rder\Gamma_{\underline{I}}
    \Rder\Gamma_{\underline{J}}\] is in general false if the families
  $\underline{I}$ and $\underline{J}$ are not sufficiently independent from one
  another: For $R=\Z$, $\underline{I}=\underline{J}=(n_i\Z)_i$ any descending
  sequence of non-trivial ideals and $M=\Q/\Z$, the five-term-sequence in
  cohomology would start as follows:
  \[ 0 \ra \varinjlim_{i,j} \Ext^1_{\Z}(\Z/n_i, \Hom(\Z/n_j, \Q/\Z))
    \ra \varinjlim_i \Ext^1_{\Z}(\Z/n_i, \Q/\Z) \ra \dots \]
  But clearly $\Ext^1_{\Z}(\Z/n_i, \Q/\Z)=0$ and
  \[\Ext^1_{\Z}(\Z/n_i, \Hom(\Z/n_j, \Q/\Z)) = \Ext^1_{\Z}(\Z/n_i,
    \Z/n_j) = \Z/(n_i,n_j),\]
  hence
  \[\varinjlim_{i,j} \Ext^1_{\Z}(\Z/n_i, \Hom(\Z/n_j, \Q/\Z)) =
    \varinjlim_i \Z/n_i,\] so the sequence above cannot possibly be
  exact.

  This argument of course generalizes: Were
  $\Rder\Gamma_{\underline{I}+\underline{J}} =
  \Rder\Gamma_{\underline{I}} \Rder\Gamma_{\underline{J}}$, then
  \cite[chapter XV, theorem 5.12]{MR0077480} implied that
  \[ \varinjlim_i \Ext^p_R(R/I_i, \varinjlim_j \Hom_R(R/J_j, Q)) = 0\] for all
  $p>0$ and $Q$ injective, i.\,e., if the isomorphism in the derived category
  holds, then because $\Rder\Gamma_{\underline{J}}$ mapped injective objects to
  $\Gamma_{\underline{I}}$-acyclics. Using
  \cref{lm:tor-m-n-zero-then-ext-m-hom-n-q-zero}, a sufficient criterion for
  that to happen is that the transition maps eventually factor through
  $\Ext^p_R(R/\widetilde{I}_i, \Hom_R(R/\widetilde{J}_j, Q))$ for some
  $\widetilde{I}_i, \widetilde{J}_j$ with $\Tor_p^R(R/\widetilde{I}_i,
  R/\widetilde{J}_j)=0$ for all $p>0$ and this criterion appears to be close to
  optimal.
  The following proposition is a simple application of this
  principle.
\end{remark}

\begin{proposition}\label{prop:spectral-sequence-for-local-cohomology-and-sum-of-ideals}
  Let $R$ be a commutative ring and
  $(r_1,\dots,r_k,s_1,\dots,s_l)$ such a regular sequence, that
  $(s_1,\dots,s_l)$ is itself again regular. Then for the ideals
  $I=(r_i)_i$ and $J=(s_i)_i$ we have
  \[\Rder\Gamma_{\underline{I}+\underline{J}}=\Rder\Gamma_{\underline{I}}\Rder\Gamma_{\underline{J}}
    = \Rder\Gamma_{\underline{J}}\Rder\Gamma_{\underline{I}} .\]
  \begin{proof}
    Denote the ideal generated by $(r_i^m)_i$ by $I^{(m)}$ and
    similarly for $J$. Then the transition maps in the
    system
    \[\varinjlim_i \Ext^p_R(R/I^i, \varinjlim_j \Hom_R(R/J^j, Q)) = \varinjlim_{i,j}\Ext^p_R(R/I^i, \Hom_R(R/J^j, Q))\]
    eventually factor through $\Ext^p_R(R/I^{(n)}, \Hom_R(R/J^{(n)}, Q))$. But
    as we saw before, the sequences $(r_1^m,\dots,r_k^m,s_1^m,\dots,s_l^m)$ and
    $(s_1^m,\dots,s_l^m)$ are again regular, so this vanishes by
    \cref{lm:tor-m-n-zero-then-ext-m-hom-n-q-zero,prop:tor-vanishes-regular-sequences}
    for $p\geq 1$. As $\Tor^R_\bullet(-,-)$ is symmetrical for commutative
    rings, the same argument also applies for
    $\Rder\Gamma_{\underline{J}}\Rder\Gamma_{\underline{I}}.$
  \end{proof}

\end{proposition}

\section{(Avoiding) Matlis Duality}
\label{sec:matlis-duality}

First recall Pontryagin duality.

\begin{theorem}[Pontryagin duality, e.\,g.{\cite[(1.1.11)]{MR2392026}}]
  The functor $\Pi=\Hom_{\mathrm{cts}}(-,\R/\Z)$ induces a contravariant
  auto-equivalence on the category of locally compact Hausdorff abelian groups
  and interchanges compact with discrete groups. The isomorphism $A\ra
  \Pi(\Pi(A))$ is given by $a \mapsto (\varphi\mapsto \varphi(a))$.

  If $G$ is pro-$p$, then
  $\Pi(G)=\Hom_{\mathrm{cts}}(G,\Q_p/\Z_p)$. If $D$ is a discrete
  torsion group or a topologically finitely generated profinite group,
  then $\Pi(D)=\Hom_{\Z}(D,\Q/\Z)$.

  We will write $-^\vee$ for $\Pi$ if it is notationally more convenient.
\end{theorem}

Matlis duality is commonly stated as follows:

\begin{theorem}[Matlis duality, {\cite[theorem 3.2.13]{MR1251956}}]\label{thm:matlis-duality}
  Let $R$ be a complete Noetherian commutative local ring with maximal
  ideal $\mathfrak m$ and $\mathcal E$ a fixed injective hull of the
  $R$-module $R/\mathfrak m$.  Then $\Hom_R(-,\mathcal E)$ induces an
  equivalence between the finitely generated modules and the Artinian
  modules with inverse $\Hom_R(-,\mathcal E)$.
\end{theorem}

\begin{example}
  If $R$ is a discrete valuation ring, then $Q(R)/R$ is an injective
  hull of its residue field.
\end{example}

Matlis duality -- using an \emph{abstract} dualizing module instead of
a topological one~-- behaves very nicely in relation to local
cohomology. In applications however the Matlis module $\mathcal E$ is
cumbersome and in general not particularly easy to construct.

\begin{example}
  Consider the rings $R=\Z_p$, $S_1=\Z_p[\pi]$ and $S_2=\Z_p[[T]]$
  with $\pi$ a uniformizer of $\Q_p(\sqrt{p})$. Clearly the
  homomorphisms $R\ra S_i$ are local and flat and their respective
  residue fields agree.  But while
  $\mathcal E_R=\Q_p/\Z_p, \mathcal E_{S_1}\cong \Q_p/\Z_p^{\oplus 2}$
  as an abelian group. Furthermore,
  $\mathcal E_{S_2}\cong \bigoplus_{\N} \Q_p/\Z_p$ as an abelian
  group by \cref{prop:matlis-module-is-pont-of-ring}.
\end{example}

The best we can hope for in general is the following.

\begin{proposition}[{\cite[Tag 08Z5]{stacks-project}}]
  Let $R\ra S$ be a flat and local homomorphism between Noetherian
  local rings with respective maximal ideals $\mathfrak m$ and
  $\mathfrak M$. Assume that $R/\mathfrak m^n \cong S/\mathfrak M^n$
  for all $n$. Then an injective hull of $S/\mathfrak M$ as an
  $S$-module is also an injective hull of $R/\mathfrak m$ as an
  $R$-module.
\end{proposition}

Starting with pro-$p$ local rings, Matlis modules are however
intimately connected with Pontryagin duality.

\begin{lemma}\label{prop:qp-mod-zp-is-somewhat-dualizing-for-r-mod-m}
  Let $R$ be a pro-$p$ local ring with maximal ideal $\mathfrak m$.
  Then there exists an isomorphism
  of $R$-modules
  $R/\mathfrak m\cong\Hom_{\Z_p}(R/\mathfrak m,\Q_p/\Z_p)$.
  \begin{proof}
    As $R/\mathfrak m$ is finite and hence a commutative field, both
    objects are isomorphic as abelian groups. As vector spaces of the
    same finite dimension over $R/\mathfrak m$ they are hence
    isomorphic as $R/\mathfrak m$-modules and thus as $R$-modules.
  \end{proof}
\end{lemma}

\begin{lemma}\label{prop:pont-is-dual-wrt-r-vee}
  Let $R$ be a pro-$p$ local ring with maximal ideal $\mathfrak m$ and
  $M$ a finitely presented or a discrete $R$-module. Then
  $\Pi(M)=\Hom_R(M,\Pi(R))$.
  \begin{proof}
    Let first $M=\varinjlim_i M_i$ be an arbitrary direct limit of
    finitely presented $R$-modules. Then
    by~\cref{lm:hom2-commutes-with-dir-lim-iff-m-fin-pres}
    \begin{align*}
      \Hom_R(\varinjlim_i M_i, \Pi(R)) &= \varprojlim_i
                                         \Hom_R(M_i,\varinjlim_k
                                         \Pi(R/\mathfrak
                                         M^k)) \\
                                       &= \varprojlim_i \varinjlim_k
                                         \Hom_R(M_i,\Hom_{\Z_p}(R/\mathfrak M^k,
                                         \Q_p/\Z_p))\\
                                       &= \varprojlim_i \varinjlim_k \Hom_{\Z_p}(M_i/\mathfrak M^k,
                                         \Q_p/\Z_p) \\
                                       &= \varprojlim_i \Pi(M_i).
    \end{align*}
    If $M$ itself is finitely presented, this shows the
    proposition. If $M$ is discrete, we can take the $M_i$ to be
    discrete and finitely presented (i.\,e., finite). The projective
    limit of their duals exists in the category of compact $R$-modules
    and it follows that $\varprojlim_i \Pi(M_i)=\Pi(M)$.
  \end{proof}
\end{lemma}

\begin{proposition}\label{prop:matlis-module-is-pont-of-ring}
  Let $R$ be a Noetherian pro-$p$ local ring with maximal
  ideal $\mathfrak m$. Then $\Pi(R)=\Hom_{\mathrm{cts}}(R,\Q_p/\Z_p)$
  is an injective hull of $R/\mathfrak m$ as an $R$-module.
  \begin{proof}
    $\Pi(R)$ is injective as an abstract $R$-module: By Baer's
    criterion it suffices to show that
    $\Hom_R(R,\Pi(R))\ra\Hom_R(I,\Pi(R))$ is surjective for every
    (left-)ideal $I$ of $R$.  By \cref{prop:pont-is-dual-wrt-r-vee},
    this is equivalent to the surjectivity of $\Pi(R)\ra \Pi(I)$,
    which is clear.
    
    In lieu of \cref{prop:qp-mod-zp-is-somewhat-dualizing-for-r-mod-m}
    it hence suffices to show that
    $\Hom_{\Z_p}(R/\mathfrak m,\Q_p/\Z_p)\subseteq
    \Hom_{\mathrm{cts}}(R,\Q_p/\Z_p)$ is an essential extension, so
    take $H$ to be an $R$-submodule of
    $\Hom_{\mathrm{cts}}(R,\Q_p/\Z_p)$ and $0\neq f\in H$. Then by
    continuity, $f$ descends to
    \[f\colon R/\mathfrak m^{k+1}\ra\Q_p/\Z_p\] with $k$ minimal. It follows
    that there exists an element $r\in\mathfrak m^k$ with $f(r)\neq 0$.
    $rf\colon R\ra\Q_p/\Z_p$ is consequentially also not zero, lies in $H$ but
    now descends to \[rf\colon R/\mathfrak m\ra \Q_p/\Z_p,\]i.\,e., $H\cap
    \Hom_{\Z_p}(R/\mathfrak m,\Q_p/\Z_p)\neq 0$.
  \end{proof}
\end{proposition}

\begin{corollary}\label{prop:pont-is-matlis-for-fg-and-art-mods}
  Let $R$ be a commutative pro-$p$ Noetherian commutative local
  ring. Then if $M$ is finitely generated or Artinian, Matlis and
  Pontryagin duality agree.
  \begin{proof}
    Immediate from
    \cref{prop:pont-is-dual-wrt-r-vee,prop:matlis-module-is-pont-of-ring}.
  \end{proof}
\end{corollary}


\begin{proposition}\label{prop:loc-coh-is-lim-of-pont-of-shift-of-invariants-of-pont}
  Let $R$ satisfy Matlis duality via $T=\Hom_R(-,\mathcal E)$. Let
  $\underline{I}$ be a decreasing family of ideals generated by
  regular sequences of length $d$. Then
  \[ \Rder\Gamma_{\underline{I}} = \varinjlim_n T \circ [d] \circ
    \RHom(R/ I_n, -) \circ T\] on $\Dder^+_c(\RMod{R}).$
  \begin{proof}
    By \cref{prop:rhom-r-mod-regular-commutes-with-matlis} it follows
    that
    $\Rder\Gamma_{\underline{I}} = \varinjlim_n \RHom(R/I_n, -) \circ
    T \circ T = \varinjlim_n T \circ [d]\circ\RHom(R/I_n, -)
    \circ T$.
  \end{proof}
\end{proposition}

\section{Tate Duality and Local Cohomology}


\begin{remark}
  Working in the derived category makes a number of subtleties more explicit
  than working only with cohomology groups. Assume that $R$ is a complete local
  commutative Noetherian ring with finite residue field of characteristic $p$
  and $G$ an analytic pro-$p$ group. Then every $\Lambda=R[[G]]$-module has a
  natural topology via the filtration of augmentation ideals of $\Lambda$. It is
  furthermore obvious to consider the following two categories:
  \begin{itemize}
  \item $\mathcal C(\Lambda)$, the category of compact
    $\Lambda$-modules (with continuous $G$-action),
  \item $\mathcal D(\Lambda)$, the category of discrete $\Lambda$-modules (with
    continuous $G$-action).
  \end{itemize}
  Pontryagin duality then induces equivalences between $\mathcal C(\Lambda)$ and
  $\mathcal D(\Lambda^\circ)$, where $-^\circ$ denotes the opposite ring. It is
  furthermore well-known that both categories are
  abelian,
  that $\mathcal C(\Lambda)$ has exact projective limits and enough projectives
  and analogously that $\mathcal D(\Lambda)$ has exact direct limits and enough
  injectives. It is important to note that the notion of continuous
  $\Lambda$-homomorphisms and abstract ones often coincides: If $M$ is finitely
  generated with the quotient topology and $N$ is either compact or discrete,
  every $\Lambda$-homomorphism $M\ra N$ is continuous, cf.~\cite[lemma 3.1.4]{MR2954997}.

  In what follows we want to compare Tate cohomology, i.\,e.\ $\Lder D$ as
  defined below, with other cohomology theories such as local cohomology. Now
  Tate cohomology is defined on the category of \emph{discrete} $G$-modules and
  we hence have a contravariant functor $\Lder D\colon \Dder^+(\mathcal
  D(\Lambda)) \ra \Dder^-(\RMod{\Lambda^\circ})$. Local cohomology on the other
  hand is defined on $\Dder^+(\RMod{\Lambda})$ or any subcategory that contains
  sufficiently many acyclic (e.\,g.\ injective) modules. This is not necessarily
  satisfied for $\Dder^+(\mathcal C(\Lambda))$. A statement such as
  \[\Lder D \circ \Pi = [d] \circ \Rder\Gamma_{\underline{I}}\]
  without further context hence does not make much sense: The
  implication would be that this would be an isomorphism of functors
  defined on $\Dder^b(\mathcal C(\Lambda))$, but
  $\Rder\Gamma_{\underline{I}}$ doesn't exist on
  $\Dder^b(\mathcal C(\Lambda))$.
\end{remark}

\begin{definition}
  Let $G$ be a profinite group and $A$ a discrete $G$-module. Denote
  with $D$ the functor \[D\colon A \mapsto \varinjlim_U (A^U)^*\]
  where $N^*=\Hom_{\Z}(N,\Q/\Z)$, the limit runs over the open normal
  subgroups of $G$ with the dual of the corestriction being the
  transition maps (cf.~\cite[II.5 and III.4]{MR2392026}). $D$ is right
  exact and contravariant and $D(A)$ has a continuous action of $G$
  \emph{from the right}. Denote its left derivation in the derived
  category of discrete $G$-modules by
  \[\Lder D\colon \Dder^+(\mathcal D(\widehat{\Z}[[G]])) \ra
    \Dder^-(\RMod{\widehat{\Z}[[G]]^\circ})\] (where $-^\circ$ denotes
  the opposite ring), so $D_i(A)= \Lder^{-i} D(A) = \varinjlim_U H^i(U,A)^*$.

  If $G$ is a profinite group, $R$ a profinite ring and $A$ a discrete
  $R[[G]]$-module, then $D(A)$ is again an $R[[G]]$-module, so
  \[\Lder D\colon \Dder^+(\mathcal D(R[[G]]))\ra
  \Dder^-(\RMod{R[[G]]^\circ}),\] where $\mathcal D(R[[G]])$ denotes the
  category of discrete $R[[G]]$-modules. Furthermore, we can of course
  also look at the functor
  \[\Lder D\colon \Dder^+(\RMod{R[[G]]})\ra
    \Dder^-(\RMod{R[[G]]^\circ}).\] Naturally, these functors don't
  necessarily coincide.
\end{definition}

\begin{proposition}\label{prop:inj-objects-in-discrete-modules-indep-of-coeff}
  Let $R$ be such a profinite ring, that the structure morphism
  $\widehat{\Z}\ra R$ gives it the structure of a finitely presented
  flat $\widehat{\Z}$-module. Let $G$ be a profinite group such that
  $R[[G]]$ is a Noetherian local ring with finite residue field. (This
  is the case if $G$ is a $p$-adic analytic group and $R$ is the
  valuation ring of a finite extension over $\Z_p$.)

  Then an injective discrete $R[[G]]$-module is an injective discrete
  $G$-module.
  \begin{proof}
    By
    \cref{prop:right-adj-preserves-inj-obj-if-left-adj-exact} it
    suffices to show that
    $?\colon\mathcal D(R[[G]])\ra \mathcal D(\widehat{\Z}[[G]])$ has
    an exact left adjoint. It is clear that
    \[ M \mapsto R[[G]] \otimes_{\widehat{\Z}[[G]]} M \] is an algebraic exact
    left adjoint, so it remains to show that $R \otimes_{\widehat{\Z}} M =
    R[[G]]\otimes_{\widehat{\Z}[[G]]} M$ is a discrete $R[[G]]$-module. Now $M$
    is the direct limit of finite modules, hence so is $ R
    \otimes_{\widehat{\Z}} M$. But for a finite $R[[G]]$-module $N$ this is
    clear as then $\mathfrak m^k M_i=0$ for some $k$ with $\mathfrak m$ the
    maximal ideal of $R[[G]]$ by Nakayama.
  \end{proof}
\end{proposition}

\begin{corollary}\label{prop:tate-d-independent-of-coefficients-if-fg-flat-etc}
  The following diagram commutes if $R$ is a finitely presented flat
  $\widehat{\Z}$-module with $R[[G]]$ Noetherian and local with finite
  residue field:
  \[ \begin{tikzcd}
      \Dder^+(\mathcal D(R[[G]])) \arrow{d}{\Lder D} \arrow{r}{?} &
      \Dder^+(\mathcal D(\widehat{\Z}[[G]]))
      \arrow{d}{\Lder D} \\
      \Dder^-(\RMod{R[[G]]^\circ})
      \arrow{r}{?}  &     \Dder^-(\RMod{\widehat{\Z}[[G]]^\circ}) 
    \end{tikzcd}\]
  \begin{proof}
    Clearly the forgetful functors and $D$ all commute on the level of
    categories of modules. The result then follows from
    \cref{prop:inj-objects-in-discrete-modules-indep-of-coeff}.
  \end{proof}
\end{corollary}

\begin{proposition}[{\cite[corollary 3.1.6, proposition 3.1.8]{MR2954997}}]
  \label{prop:fg-proj-is-compact-proj}

  \begin{enumerate}
  \item Let $M$ be an Artinian $\Lambda$-module. Then $\Lambda\times M\ra M$
  is continuous if we give $M$ the discrete topology.
\item Let $N$ be a Noetherian $\Lambda$-module. Then $N$ is compact if we
  give it the topology induced by $\Lambda$.
\item The functor
  \[\cat{f}.\cat{g}.\textrm{-}\RMod{\Lambda}\ra\mathcal C(\Lambda)\]
  maps projective objects to projectives.
  \end{enumerate}
\end{proposition}

\begin{proposition}\label{prop:local-coh-of-torus-is-tate-coh-of-matlis}
  For the $d$-dimensional torus $G\cong\Z_p^d$ and its Iwasawa algebra
  $\Lambda(G)=\varprojlim_i \mathcal O[G/G^{p^i}]$ over a discrete
  complete valuation ring $\mathcal O$ with residue characteristic $p$
  and uniformizer $\pi$ and with augmentation
  ideals $I_i = \ker \Lambda(G) \ra \mathcal O[G/G^{p^i}]$ the
  following holds in $\Dder^b_c(\RMod{\Lambda})$:
  \[\Lder D \circ T \cong [d] \circ \Rder\Gamma_{\underline{I}},\]
  Especially the following diagram commutes:
  \[ \begin{tikzcd} \Dder^b_c(\RMod{\Lambda}) \arrow[dash]{r}{\cong}
      \arrow[hook]{d} & \Dder^b(\cat{f}.\cat{g}.\textrm{-}\RMod{\Lambda})
      \arrow{r} & \Dder^b(\mathcal C(\Lambda)) \arrow{r}{T} &
      \Dder^b(\mathcal
      D(\Lambda^\circ)) \arrow{d}{\Lder D} \\
      \Dder^b(\RMod{\Lambda})
    \arrow{rrr}{[d]\circ\Rder\Gamma_{\underline{\mathfrak
          m}}}&&& \Dder^b(\RMod{\Lambda})
    \end{tikzcd} \]

  \begin{proof}
    $\Lambda(G)$ is a regular local ring with maximal ideal
    $(\pi, \gamma_1-1,\dots, \gamma_d-1)$ for any set of topological
    generators $(\gamma_i)_i$ of $G$. One immediately verifies that
    the sequences $\gamma_1^{p^i}-1,\dots, \gamma_d^{p^i}-1$ are again
    regular and generate the ideals $I_i$.

    By
    \cref{prop:loc-coh-is-lim-of-pont-of-shift-of-invariants-of-pont},
    \[ [d] \circ \Rder\Gamma_{\underline{I}}= \varinjlim_n
      T\circ \RHom_\Lambda(\Lambda/I_n,-) \circ T.\]

    Take a bounded complex $M$ of finitely generated $R$-modules that
    is quasi-isomorphic to a bounded complex $P$ of finitely generated
    projective modules. The resulting complex $T(P)$ is then not
    only a bounded complex of injective discrete modules by Pontryagin
    duality and \cref{prop:fg-proj-is-compact-proj}, but also a
    bounded complex of injective abstract $\Lambda$-modules by
    \cref{lm:hom-of-flat-inj-is-inj}. In all relevant derived
    categories $T(M)\cong T(P)$ holds. As
    $\Hom_\Lambda(\Lambda/I_n,-)=(-)^{G^{p^n}}$ by construction, we
    can compute $[d] \circ \Rder\Gamma_{\underline{I}} (M)$ as follows
    (keeping \cref{prop:pont-is-matlis-for-fg-and-art-mods} in mind):
    \begin{align*}
      [d] \circ \Rder\Gamma_{\underline{I}} (M) &= [d] \circ
                                              \Rder\Gamma_{\underline{I}}
                                              (P) \\
                        &= \varinjlim_n T\circ \RHom_\Lambda(\Lambda/I_n,-)
                                              \circ  T(P)\\
                        & = \varinjlim_n T(\Hom_\Lambda(\Lambda/I_n,T(P))) \\
                        & = \varinjlim_n T(T(P)^{G^{p^n}}) \\
                        & = \varinjlim_n \Pi(T(P)^{G^{p^n}}) \\
                        & = D(T(P)) = \Lder D \circ T(M).
    \end{align*}
  \end{proof}
\end{proposition}

\begin{lemma}\label{prop:rhom-r-mod-m-k-restricts-to-der-coh-lambda-mod}
  Let $R$ be a commutative Noetherian ring with unit and $R\ra S$ a flat ring
  extension with $R$ contained in the centre of $S$ and $S$ again
  (left-)Noetherian.

  Then \[\RHom_R\colon \Dder^-(\RMod{R})^\opp \times\Dder^+(\RMod{R}) \ra
    \Dder^+(\RMod{R})\]
  extends to
  \[\RHom_R\colon \Dder^-(\RMod{R})^\opp \times\Dder^+(\RMod{S}) \ra
    \Dder^+(\RMod{S}),\]
  which in turn restricts to \[\RHom_R\colon \Dder^b_c(\RMod{R})^\opp \times\Dder^b_c(\RMod{S}) \ra
    \Dder^b_c(\RMod{S}),\]
    \begin{proof}
      First note that if $M$ is an $R$-left-module and $N$ an
      $S$-left-module, then $\Hom_R(M,N)$ carries the structure of an
      $S$-left-module via $(s f)(m) = sf(m)$. Then
      $\Hom_R(R,S)\cong S$ as $S$-left-modules and the
      following diagram commutes:\[\begin{tikzcd} \RMod{S}
          \arrow{rr}{\Hom_R(M,-)}\arrow{d}{?} && \RMod{S}\arrow{d}{?} \\
          \RMod{R}\arrow{rr}{\Hom_R(M, -)} && \RMod{R}\end{tikzcd} \] As
      $S$ is a flat $R$-module, $?$ preserves injectives by
      \cref{prop:right-adj-preserves-inj-obj-if-left-adj-exact} and we can
      compute $\RHom_R(M,-)$ in either category.

      If $M$ is a finitely generated $R$-module and $N$ a finitely generated
      $S$-module, then $\Hom_R(M,N)$ is again a finitely generated $S$-module,
      as $S$ is left-Noetherian. If $M$ is a bounded complex of finitely
      generated $R$-modules, then it is quasi-isomorphic to a bounded complex of
      finitely generated projective $R$-modules. The result then follows at
      once.
    \end{proof}
\end{lemma}

\begin{remark}
  Note however that $\RHom_R$ does \emph{not} extend to a
  functor \[\RHom_R\colon \Dder^-(\RMod{S})^\opp \times\Dder^+(\RMod{R}) \ra
    \Dder^+(\RMod{S}).\] Even in those cases where we can give $\Hom_R(M,A)$ the
  structure of an $S$-module (e.\,g.\ when $S$ has a Hopf structure with
  antipode $s\mapsto\overline{s}$ via $(sf)(m)=f(\overline{s}m)$), projective
  $S$-modules in general are not projective. This is specially true for
  $R[[G]]$, which is a flat, but generally not a projective $R$-module. 
\end{remark}


An essential ingredient in the proof of this section's main theorem is
Grothendieck local duality. It is commonly stated as follows:

\begin{theorem}[Local duality, {\cite[V.6.2, V.9.1]{MR0222093}}]
  \label{thm:local-duality}
  Let $R$ be a commutative regular local ring of dimension $d$ with
  maximal ideal $\mathfrak m$, and $\mathcal E$ a fixed injective hull
  of the $R$-module $R/\mathfrak m$. Denote with $R[d]$ the complex
  concentrated in degree $-d$ with entry $R$. Then
  \[ \Rder \Gamma_{\underline{\mathfrak m}} \cong T \circ \RHom_R(-,
    R[d]) = [-d] \circ T \circ \RHom_R(-,R) \] on $\Dder^b_c(\RMod{R})$.
\end{theorem}

The regularity assumption on $R$ can be weakened if one is willing to
deal with a dualizing complex that is not concentrated in just one
degree (loc.\,cit.). Relaxing commutativity however is more subtle and
will be the focus of \cref{sec:local-dual-iwas}.

\begin{theorem}\label{prop:jannsen-spec-seq-in-der-cat}
  Let $\mathcal O$ be a pro-$p$ discrete valuation ring,
  $R=\mathcal O[[X_1,\dots,X_t]]$ with maximal ideal $\mathfrak m$,
  $G=\Z_p^s$ and $\Lambda = R[[G]]$. Then
  \[ T \circ \RHom_\Lambda(-,\Lambda) \cong [t+1] \circ \Rder
    \Gamma_{\mathfrak m} \circ \Lder D \circ T \] on
  $\Dder^b_c(\RMod{\Lambda})$. The right hand side can furthermore be
  expressed as
  \[ \Rder \Gamma_{\mathfrak m} \circ \Lder D \circ T \cong
    \varinjlim_k\, \Lder D \circ T \circ \RHom_R(R/\mathfrak m^k, -).\]
  \begin{proof}
    $\Lambda$ is again a regular local ring, now of dimension
    $t+s+1$. Denote its maximal ideal by $\mathfrak M$.  By
    \cref{thm:local-duality}
    \[ \Rder\Gamma_{\underline{\mathfrak M}} \cong T\circ
      \RHom_\Lambda(-,\Lambda[s+t+1]) = [-s-t-1] \circ T \circ
      \RHom_\Lambda(-,\Lambda) .\]
    Now $\mathfrak M= \mathfrak m + (\gamma_1-1,\dots,\gamma_s-1)$ and
    if $x_1,\dots,x_{t+1}$ is a regular sequence in $R$, then
    $x_1,\dots,x_{t+1}, \gamma_1-1,\dots,\gamma_s-1$ is a regular
    sequence in $\Lambda$. Furthermore, the sequence
    $\gamma_1-1,\dots,\gamma_s-1$ is of course itself regular in
    $\Lambda$. Let it generate the ideal $I$. Then we can apply
    \cref{prop:spectral-sequence-for-local-cohomology-and-sum-of-ideals},
    i.\,e.,
    \[ \Rder\Gamma_{\mathfrak M} \cong
      \Rder\Gamma_{\underline{\mathfrak m}} \circ
      \Rder\Gamma_{\underline{I}}.\] By
    \cref{prop:local-coh-of-torus-is-tate-coh-of-matlis}, we have
    $\Rder\Gamma_{\underline{I}} = [-s] \circ \Lder D \circ T$, which
    shows the first isomorphism.

    Consider furthermore the functor
    $\varinjlim_k \Lder D \circ T \circ \RHom_R(R/\mathfrak m^k,
    -)$. By \cref{prop:rhom-r-mod-m-k-restricts-to-der-coh-lambda-mod}
    we can compute it on $\Dder^b_c(\RMod{\Lambda})$ as
    \begin{align*}\varinjlim_k \Lder D \circ T \circ \RHom_R(R/\mathfrak m^k, -)
      & \cong \varinjlim_k [s] \circ\Rder \Gamma_{\underline{I}} \circ 
        \RHom_R(R/\mathfrak m^k, -)\\
      &\cong [s] \circ \Rder \Gamma_{\underline{I}} \circ \varinjlim_k
        \RHom_R(R/\mathfrak m^k, -) \\
      & = [s] \circ \Rder \Gamma_{\underline{I}} \circ
        \Rder\Gamma_{\underline{\mathfrak m}} \\
      &\cong [s] \circ \Rder \Gamma_{\underline{\mathfrak m}} \circ
        \Rder\Gamma_{\underline{I}}\\
      &\cong [s] \circ \Rder \Gamma_{\underline{\mathfrak m}} \circ
        [-s] \circ \Lder D \circ T = \Rder \Gamma_{\underline{\mathfrak m}}\circ \Lder D \circ T, \end{align*}
      as by \cref{prop:local-coh-commutes-with-direct-limits}, local
      cohomology commutes with direct limits, $\Rder\Gamma_{\mathfrak m}$ and
      $\Rder\Gamma_{\underline{I}}$ commute by \cref{prop:spectral-sequence-for-local-cohomology-and-sum-of-ideals} and two choice
      applications of 
      \cref{prop:local-coh-of-torus-is-tate-coh-of-matlis}.
  \end{proof}
\end{theorem}

\begin{remark} 
  If we express \cref{prop:jannsen-spec-seq-in-der-cat} in terms of a
  spectral sequence, it looks like this:
  \[ \varinjlim_k\Lder^pD(T(\Ext^q_R(R/\mathfrak m^k, M))) \Longrightarrow
    T(\Ext^{t+1-p-q}_\Lambda(M,\Lambda)).\] Writing $\E^\bullet_\Lambda$ for
  $\Ext^\bullet_\Lambda(-,\Lambda)$, flipping the sign of $p$ and shifting
  $q\mapsto t+1-q$ then yields
  \[ \varinjlim_k D_p( \Ext^{t+1-q}_R(R/\mathfrak m^k, M)^\vee)
    \Longrightarrow \E^{p+q}_\Lambda(M)^\vee \] and the following exact five
  term sequence:

  \[ \begin{tikzcd}\E^2_\Lambda(M)^\vee \arrow{r} & \varinjlim_k D_2(\Ext_R^{t+1}(R/\mathfrak m^k,
    M)^\vee) \arrow{r} & \varinjlim_k D(\Ext_R^{t}(R/\mathfrak m^k,
    M)^\vee)  \arrow[out=355, in=175, looseness=1.5, overlay]{lld} \\ \E^1_\Lambda(M)^\vee \arrow{r} & \varinjlim_k D_1(\Ext_R^{t+1}(R/\mathfrak m^k,
    M)^\vee) \arrow{r} & 0\end{tikzcd}\]

\end{remark}

\section{Iwasawa Adjoints}
\label{sec:iwasawa-adjoints}

In this section let $R$ be a pro-$p$ commutative local ring with maximal ideal
$\mathfrak m$ and residue field of characteristic $p$. Let $G$ be a compact
$p$-adic Lie group and $\Lambda=\Lambda(G)=\varprojlim_U R[[G/U]]$, where $U$
ranges over the open normal subgroups of $G$. As is customary, we set again
$\E^\bullet_\Lambda(M) = \Ext^\bullet_\Lambda(M,\Lambda)$.

Even though $R$ is commutative, we will differentiate between $R$- and
$R^\circ$-modules, as we will regard $R$ as a subring of $\Lambda$.

\begin{remark}
  Note that if $M$ is a left $\Lambda$-module, it also has an operation of
  $\Lambda$ from the right given by $mg=g^{-1}m$. This of course does not give
  $M$ the structure of a $\Lambda$-bimodule, as the actions are not compatible.
  We can however still give $\Hom_\Lambda(M,\Lambda)$ the structure of a left
  $\Lambda$-module by $(g.\varphi)(m)= \varphi(m)g^{-1}$.
\end{remark}

\begin{lemma}\label{prop:e0-as-proj-lim-over-coinv}
  $\E^0_\Lambda(M) = \varprojlim_U \Hom_R(M_U, R)$ for finitely
  generated $\Lambda$-modules $M$, where the transition map for a pair
  of open normal subgroups $U\leq V$ are given by the dual of the
  trace map $M_V\ra M_U, m\mapsto\sum_{g\in V/U} gm$.
  \begin{proof}
    Note first that as $\Hom_\Lambda(M,-)$ commutes with projective
    limits,
    \[\Hom_\Lambda(M,\Lambda) = \varprojlim_U \Hom_\Lambda(M, R[G/U]) =
      \varprojlim_U \Hom_{R[G/U]}(M_U, R[G/U]).\] For $U$ an open
    normal subgroups of $G$, consider the trace map
    \begin{gather*}
      \Hom_R(M_U, R)\ra\Hom_{R[G/U]}(M_U, R[G/U])\\
      \varphi\mapsto \left(m\mapsto\sum_{g\in G/U}  
        \varphi(g^{-1}m)\cdot g\right),\end{gather*} which is clearly
    an isomorphism of $R$-modules and induces the required isomorphism
    to the projective system mentioned in the proposition.
  \end{proof}
\end{lemma}

\begin{proposition}\label{prop:pont-of-iw-adjoints-is-colim-pont-rhom-inv}
  On $\Dder^b_c(\RMod{\Lambda})$ we have
  \[\Pi \circ
    \RHom_\Lambda(-, \Lambda) \cong \varinjlim_U \Pi \circ
    \RHom_R(-,R) \circ \Lder(-)_{U}.\]
  \begin{proof}
    Immediate by \cref{prop:e0-as-proj-lim-over-coinv}, as $(-)_U$
    clearly maps finitely generated free $\Lambda$-modules to finitely
    generated free $R$-modules.
  \end{proof}
\end{proposition}

\begin{remark}
  The spectral sequence attached to
  \cref{prop:pont-of-iw-adjoints-is-colim-pont-rhom-inv} looks like this: \[
    \varinjlim_U \Ext_R^p(H_{q}(U, M), R)^\vee \Longrightarrow
    \E^{p+q}_\Lambda(M)^\vee\]
  Its five term exact sequence is given by

  \[ \begin{tikzcd}
      \E^2_\Lambda(M)^\vee \arrow{r} & \varinjlim_U \Ext^2_R(M_U, R)^\vee
      \arrow{r} & \varinjlim_U \Hom_R(H_1(U,M), R)^\vee \arrow[looseness=1.5,
      overlay, out=355, in=175]{dll} \\
      \E^1_\Lambda(M)^\vee \arrow{r} & \varinjlim_U \Ext^1_R(M_U, R)^\vee
      \arrow{r} & 0
    \end{tikzcd}\]
\end{remark}
\begin{lemma}\label{prop:iw-dual-as-colim-of-coinv-otimes-r-vee}
  $\Hom_\Lambda(M,\Lambda)^\vee \cong \varinjlim_U R^\vee \otimes_R M_U$
  for finitely generated $\Lambda$-modules $M$.
  \begin{proof}
    \begin{align*}
      \Hom_\Lambda(M,\Lambda)^\vee & \cong \varinjlim_U \Pi(\Hom_R(M_U, R)) \\
                                 & \cong\varinjlim_U \Pi\Hom_{R^\circ}(\Pi(R),
                                   \Pi(M_U))\\
                                 &\cong \Pi \circ \Pi(  \varinjlim_U M_U\otimes_{R^\circ}
                                   \Pi(R))\\
                                 &\cong \varinjlim M_U \otimes_{R^\circ}
                                   R^\vee.
    \end{align*}
  \end{proof}
\end{lemma}

\begin{proposition}
  \label{prop:pont-of-iw-adjoints-is-tor-tate-pont}
  $\Pi \circ \RHom_\Lambda(-,\Lambda) \cong \left(R^\vee \otimes^{\cat{L}}_R
  -\right) \circ \Lder D \circ \Pi$ on $\Dder^b_c(\RMod{\Lambda})$.
  \begin{proof}
    Using
    \cref{prop:e0-as-proj-lim-over-coinv,prop:iw-dual-as-colim-of-coinv-otimes-r-vee},
    usual Pontryagin duality (where we occasionally write $-^\vee$ for
    the dual), and the fact that tensor products commute with direct
    limits, we get:
    \begin{align*}
      \Hom_\Lambda(M,\Lambda)^\vee & \cong  \varinjlim_U R^\vee \otimes_R
                                     M_U\\
                                   &\cong  R^\vee \otimes_R \varinjlim \Pi(\Pi(M)^U)\\
                                   &= R^\vee \otimes_R D(\Pi(M)).
    \end{align*}
    It hence remains to show that $(D\circ \Pi)$ maps projective
    objects to $R^\vee \otimes_R -$-acyclics and it actually suffices
    to check this for the module $\Lambda$. But
    $D(\Pi(\Lambda))=\varinjlim_U R[G/U]$ is clearly
    $R^\vee \otimes_R -$-acyclic.
  \end{proof}
\end{proposition}

    

\begin{remark}
  The spectral sequence attached to
  \cref{prop:pont-of-iw-adjoints-is-tor-tate-pont} looks like this:
  \[ \Tor_{p}^R(R^\vee,D_{q}(M^\vee)) \Longrightarrow
    \E^{p+q}_\Lambda(M)^\vee,\] which yields the following five term exact
  sequence:

  \[\begin{tikzcd}
      \E^2_\Lambda(M)^\vee \arrow{r} & \Tor_2^R(R^\vee, D(M^\vee)) \arrow{r} &
      R^\vee \otimes_R D_1(M^\vee) \arrow[looseness=1.5, overlay, out=355, in=175]{lld} \\
      \E^1_\Lambda(M)^\vee \arrow{r} & \Tor_1^R(R^\vee, D(M^\vee)) \arrow{r} & 0
    \end{tikzcd}\] This also proves that $\E^p_\Lambda(M)=0$ if $p>\dim G + \dim
  R$. If $\dim R=1$, the spectral sequence degenerates and we can compute
  $\E^p_\Lambda(\tor_R M)$ and $\E^p_\Lambda(M/\tor_R M)$ akin to
  \cite[(5.4.13)]{MR2392026}.
\end{remark}


\begin{lemma}\label{prop:r-vee-is-colim-over-residues}
  $R^\vee \cong \varinjlim_k R/\mathfrak m^k$ if $R$ is regular.
  \begin{proof}
    $R$ satisfies local duality by assumption, hence by
    \cref{prop:pont-is-matlis-for-fg-and-art-mods,prop:duality-ext-tor-for-reg-seq-classical}
    $R^\vee = T(R) \cong \Rder^d\Gamma_{\underline{\mathfrak m}}(R) =
    \varinjlim_k \Ext^d_R(R/\mathfrak m^{(k)}, R) \cong \varinjlim_k
    R/\mathfrak m^{(k)},$ where we again denote with $\mathfrak m^{(k)}$ the
  ideal generated by $k$th powers of generators of $\mathfrak m$.
  \end{proof}
\end{lemma}

\begin{lemma}\label{prop:iw-dual-as-colim-of-coniv-of-residues}
  $\Hom_\Lambda(M,\Lambda)^\vee \cong \varinjlim_{U,k} (M/\mathfrak
  m^k)_U$ for finitely generated $\Lambda$-modules $M$ and regular $R$.
  \begin{proof}
    \begin{align*}
      \Hom_\Lambda(M,\Lambda)^\vee &\cong R^\vee \otimes_R \varinjlim_U M_U \\
                                   &\cong \varinjlim_U \varinjlim_k
                                     R/\mathfrak m^k \otimes_R M_U  \\
                                   &\cong \varinjlim_U \varinjlim_k
                                     (M/\mathfrak m^k
                                     M)_U\end{align*}
                                   using \cref{prop:iw-dual-as-colim-of-coinv-otimes-r-vee,prop:r-vee-is-colim-over-residues}.
  \end{proof}
\end{lemma}

\begin{proposition}\label{prop:pont-of-iw-adjoints-is-colim-of-tate-of-pont-of-loc-coh}
  If $R$ is regular, $ \Pi \circ \RHom_\Lambda(-, \Lambda) \cong
  \varinjlim_k \Lder D \circ \Pi \circ \left(R/\mathfrak m^k \otimes^{\cat{L}}_R
    -\right) \cong \varinjlim_k \Lder D \circ \Pi \circ [d] \circ
  \RHom_R(R/\mathfrak m^k, -)$ on $\Dder^b_c(\RMod{\Lambda}).$
  \begin{proof}
    By \cref{prop:iw-dual-as-colim-of-coniv-of-residues}
    \begin{align*}
      \Hom_\Lambda(M,\Lambda)^\vee &\cong \varinjlim_U \varinjlim_k (M/\mathfrak m^k M)_U \\
                                   &\cong \varinjlim_{U,k} \Pi(\Pi(M/\mathfrak m^k M)^U)\\
      &= \varinjlim_k D(\Pi(R/\mathfrak m^k \otimes_R M)).
    \end{align*}
    By \cref{prop:duality-rhom-tor-for-reg-seq-in-der-cat} it suffices
    to show that $(\Lambda/\mathfrak m^k)^\vee$ is $D$-acyclic. But
    \[\Lder^{-i}D((\Lambda/\mathfrak m^k)^\vee) = \varinjlim_U H^i(U,
      R/\mathfrak m^k[[G]]^\vee)^* = \varinjlim_U H_i(U, R/\mathfrak
      m^k[[G]])=0\]for $i>0$ by Shapiro's Lemma.

    The other isomorphism now follows from
    \cref{prop:duality-rhom-tor-for-reg-seq-in-der-cat}.
  \end{proof}
\end{proposition}

\begin{remark}
  The spectral sequences attached to
  \cref{prop:pont-of-iw-adjoints-is-colim-of-tate-of-pont-of-loc-coh} look like
  this: \[ \varinjlim_k D_{p}(\Tor_q^R(R/\mathfrak m^k, M)^\vee)
    \Longrightarrow \E^{p+q}_\Lambda(M)^\vee\] and \[ \varinjlim_k
    D_{p}(\Ext_R^{d-q}(R/\mathfrak m^k, M)^\vee) \Longrightarrow
    \E^{p+q}_\Lambda(M)^\vee\]
  with exact five term sequences
  \[\begin{tikzcd}
      \E^2_\Lambda(M)^\vee \arrow{r} & \varinjlim_k D_2( (M/\mathfrak m^kM
      )^\vee)\arrow{r} &
      \varinjlim_k D(\Tor_1^R(R/\mathfrak m^k, M)^\vee) \arrow[looseness=1.5, overlay, out=355, in=175]{lld} \\
      \E^1_\Lambda(M)^\vee \arrow{r} & \varinjlim_k D_1((M/\mathfrak m^kM)^\vee) \arrow{r} & 0
    \end{tikzcd}\] and \[\begin{tikzcd} \E^2_\Lambda(M)^\vee \arrow{r} &
      \varinjlim_k D_2(\Ext^d_R(R/\mathfrak m^k, M)^\vee)\arrow{r} &
      \varinjlim_k D(\Ext^{d-1}_R(R/\mathfrak m^k, M)^\vee) \arrow[looseness=1.5, overlay, out=355, in=175]{lld} \\
      \E^1_\Lambda(M)^\vee \arrow{r} & \varinjlim_k D_1(\Ext^d_R(R/\mathfrak m^k
      , M)^\vee) \arrow{r} & 0
    \end{tikzcd}\]  respectively.
\end{remark}

\begin{lemma}[{\cite[proposition 3.1.7]{MR2954997}}]
  \label{prop:fg-lambda-mods-are-complete-and-sep}
  Let $M$ be a finitely generated $\Lambda$-module. Then
  $M\cong\varprojlim_k M/\mathfrak M^kM$ algebraically and
  topologically.
\end{lemma}

\begin{lemma}\label{prop:pont-of-rhom-r-mod-m-is-disc-p-tor}
  $\Pi \circ \RHom_R(R/\mathfrak m^k, -)$ maps bounded complexes of
  $\Lambda$-modules with coherent cohomology to bounded complexes whose
  cohomology modules are discrete $p$-torsion $G$-modules. If $M$ is a complex
  in $\varinjlim_U \Dder^b_c(\RMod{R[G/U]})$, then all cohomology groups of
  $\RHom_R(R/\mathfrak m^k, M)^\vee $ are furthermore finite.
  \begin{proof}
    The groups $\Ext^q_R(R/\mathfrak m^k, M)$ for $M$ finitely
    generated over $\Lambda$ are clearly $p$-torsion and finitely
    generated as $\Lambda$-modules, hence compact by
    \cref{prop:fg-lambda-mods-are-complete-and-sep}, and
    consequentially topologically profinite and pro-$p$. Their
    Pontryagin duals are thus discrete $p$-torsion $G$-modules.

    If $M$ is finitely generated over some $R[G/U]$, it is also finitely
    generated over $R$ and $\Ext^q_R(R/\mathfrak m^k, M)$ finitely generated
    over $R/\mathfrak m^k$, hence finite.
  \end{proof}
\end{lemma}

\begin{proposition} \label{prop:pont-of-iw-adjoint-with-duality-grp-has-deg-spec-seq}
  Assume that $R$ is regular. Let $G$ be a duality group
  (cf.~\cite[(3.4.6)]{MR2392026}) of dimension $s$ at $p$. Then
  \begin{align*}
    \E^m_\Lambda(M)^\vee &\cong \varinjlim_k \Lder^{-s} D
                           \Ext^{d-(m-s)}_R(R/\mathfrak m^k, M)^\vee \\
                         & \cong \varinjlim_k \Lder^{-s} D
                           \Tor_{m-s}^R(R/\mathfrak m^k,
                           M)^\vee.\end{align*}
                         for finitely
                         generated $R[G/U]$-modules $M$.
                         Especially $\Pi\circ \E^s_\Lambda$ is
                         then right-exact.
  \begin{proof}
    As $G$ is a duality group of dimension $s$ at $p$, the complex $\Lder D(M')$
    has trivial cohomology outside of degree $-s$ if $M'$ is a finite discrete
    $p$-torsion $G$-module. Together with
    \cref{prop:pont-of-rhom-r-mod-m-is-disc-p-tor} this implies that the
    spectral sequence attached to
    \cref{prop:pont-of-iw-adjoints-is-colim-of-tate-of-pont-of-loc-coh}
    degenerates and gives
    \[ \E^m_\Lambda(M)^\vee \cong \varinjlim_k \Lder^{-s} D
      \Ext^{d-(m-s)}_R(R/\mathfrak m^k, M)^\vee.\] The other isomorphism follows
    with exactly the same argument.
    Note furthermore that as $\dim R=d$, $\Ext^{d-(m-s)}_R(R/\mathfrak
    m^k, M)=0$ if $m-s<0$, hence $\E^m_\Lambda(M)^\vee=0$ if $m<s$.
  \end{proof}
\end{proposition}

\begin{theorem}\label{prop:dual-of-iw-adj-is-twist-of-local-coh-of-coeff}
  Assume that $R$ is regular and that $G$ is a Poincaré group at $p$ of
  dimension $s$ with dualizing character $\chi\colon G\ra \Z_p^\times$
  (cf.\cite[(3.7.1)]{MR2392026}), which gives rise to the “twisting functor”
  $\cat{\chi}\colon M\mapsto M(\chi)$. Assume that $R$ is a commutative complete
  Noetherian ring of global dimension $d$ with maximal ideal $\mathfrak m$. Then
  \[ T \circ \RHom_\Lambda(-,\Lambda) = \cat{\chi} \circ [d+s] \circ
    \Rder\Gamma_{\underline{\mathfrak m}}\] on
  $\varinjlim_U \Dder^b_c(\RMod{R[G/U]}).$
  \begin{proof}
    Let $I=\Q_p/\Z_p(\chi)=\cat{\chi}(\Q_p/\Z_p)$ be the dualizing
    module of $G$.  For any $p$-torsion $G$-module $A$ we have by
    \cite[(3.7)]{MR2392026} that
    $\Lder^{-s}D(A) = \varinjlim_U H^s(U,A)^* \cong \varinjlim_U
    \Hom_{\Z_p}(A,I)^U = \Hom_{\Z_p}(A,I)$, as $I$ is also a dualizing
    module for every open subgroup of $G$.
    
    Note that
    \[H^0 ( \cat{\chi} \circ [d] \circ \Rder\Gamma_{\underline{\mathfrak
          m}}) = \cat{\chi}\circ \Rder^{d}\Gamma_{\underline{\mathfrak
          m}}\] and that
    $\cat{\chi}\circ\Rder^{d}\Gamma_{\underline{\mathfrak m}}$ is hence
    right-exact. Note furthermore that
    \begin{align*}
      H^0( [-s]\circ T \circ \RHom_\Lambda(-,\Lambda)) &= \E^s(-)^\vee
      \\
                                                       &\cong \varinjlim_k \Lder^{-s}D \Ext^d_R(R/\mathfrak m^k, -)^\vee \\
                                                       &\cong \varinjlim_k \Hom_{\Z_p} (\Ext^d_R(R/\mathfrak m^k, -)^\vee,
                                                         I)\\
                                                       &\cong \cat{\chi}\circ\Rder^d\Gamma_{\underline{\mathfrak m}}
    \end{align*}
    using
    \cref{prop:pont-of-iw-adjoint-with-duality-grp-has-deg-spec-seq}
    and Pontryagin duality.

    By \cite[I.7.4]{MR0222093} the left derivation of
    $\Rder^d\Gamma_{\underline{\mathfrak m}}$ is
    $[d]\circ\Rder\Gamma_{\underline{\mathfrak m}}$: The complex
    $\Rder\Gamma_{\underline{\mathfrak m}}(R)$ is concentrated in
    degree $d$ and hence every module is a quotient of a module with
    this property, as local cohomology commutes with arbitrary direct
    limits.
  \end{proof}
\end{theorem}

Note that even though this theorem suspiciously looks like local
duality, the local cohomology on the right hand side is with respect
to the maximal ideal of the coefficient ring, not the whole Iwasawa
algebra. The local duality result is subject of the next section.

\begin{corollary}\label{prop:explicit-iw-adj}
  In the setup of \cref{prop:dual-of-iw-adj-is-twist-of-local-coh-of-coeff}
  assume that $M$ is a finitely generated $R[G/U]$-module. Then the following hold:
  \begin{enumerate}
  \item If $M$ is free over $R$, then \[ \E_\Lambda^q(M)^\vee
      \cong \begin{cases}M\otimes_R R^\vee(\chi) & \text{if } q=s\\
        0 &\text{else}\end{cases}\]
  \item If $M$ is $R$-torsion, then $\E_\Lambda^q(M)=0$ for all $q\leq s$.
    \item If $M$ is finite, then \[ \E_\Lambda^q(M)^\vee
      \cong \begin{cases}M(\chi) & \text{if } q=d+s\\
        0 &\text{else}\end{cases}\]
  \end{enumerate}
  \begin{proof}
    By \cref{prop:dual-of-iw-adj-is-twist-of-local-coh-of-coeff}, we have \[
      \E_\Lambda^q(M)^\vee \cong \Rder^{d+s-q}\Gamma_{\underline{\mathfrak
          m}}(M)(\chi)\] in any case.

    In the first case, this is just $M\otimes_R R^\vee$ for $q=s$ and zero else.
    In the second case, local duality yields
    $\Rder^d\Gamma_{\underline{\mathfrak m}}(M)\cong \Hom_R(M,R)^\vee=0$. In the
    third case, we note that $M$ has an injective resolution by modules that are
    the direct limit of finite modules.
    \Cref{prop:local-coh-commutes-with-direct-limits} together with
    \cref{prop:loc-coh-is-max-fin-submod} then imply the result.
  \end{proof}
\end{corollary}

\section{Local Duality for Iwasawa Algebras}
\label{sec:local-dual-iwas}

This section gives a streamlined proof of a local duality result for Iwasawa
algebras as first published in \citenospec{MR1924402} and generalizes the result
to more general coefficient rings. Let $G$ be a pro-$p$ Poincaré group of
dimension $s$ with dualizing character $\chi\colon G\ra\Z_p^\times$ and $R$ a
commutative Noetherian pro-$p$ regular local ring with maximal ideal $\mathfrak
m$ of global dimension $d$. Set $\Lambda=R[[G]]$.

\begin{proposition}\label{prop:g-poincare-then-top-loc-coh-is-pont-dual-for-ring}
  $\Rder\Gamma_{\underline{\mathfrak
      M}}(\Lambda)\cong\Lambda^\vee[-d-s]$ and
  $\Ext^i_\Lambda(\Lambda/\mathfrak M^l,\Lambda)\cong
  \Rder\Gamma_{\underline{\mathfrak m}}^{d+s-i}(\Lambda/\mathfrak
  M^l)(\chi).$
  
  \begin{proof}
    By
    \cref{prop:pont-of-iw-adjoint-with-duality-grp-has-deg-spec-seq},
    \[\Rder^{i}\Gamma_{\underline{\mathfrak M}}(\Lambda) =
      \varinjlim_l \E^i(\Lambda/\mathfrak M^l) \cong \varinjlim_{l}
      \left( \varinjlim_k \Lder^{-s}D (\Ext^{d-(i-s)}_R(R/\mathfrak m^k,
        \Lambda/\mathfrak M^l)^\vee)\right)^\vee.\] As in the proof of
    \cref{prop:dual-of-iw-adj-is-twist-of-local-coh-of-coeff},
    \begin{align*} \Lder^{-s}D (\Ext^{d-(i-s)}_R(R/\mathfrak m^k,
      \Lambda/\mathfrak M^l)^\vee) &\cong
      \Hom_{\Z_p}(\Ext^{d-(i-s)}_R(R/\mathfrak m^k, \Lambda/\mathfrak
      M^l)^\vee, \Q_p/\Z_p(\chi))\\ &\cong
      \Ext^{d-(i-s)}_R(R/\mathfrak m^k, \Lambda/\mathfrak
      M^l)(\chi)\end{align*} and in the direct limit over $k$ this
    becomes
    $\Rder\Gamma^{d-(i-s)}_{\underline{\mathfrak m}}(\Lambda/\mathfrak
    M^l)(\chi)$.

    $\Gamma_{\underline{\mathfrak m}}$ restricted to the subcategory
    of finite $\Lambda$- (or $R$-)modules is the identity. As
    $\Gamma_{\underline{\mathfrak m}}$ commutes with arbitrary direct
    limits, this is also true for the category of discrete
    $\Lambda$-modules. As the latter category contains sufficiently
    many injectives,
    $\Rder\Gamma_{\underline{\mathfrak
        m}}(N)=N$ if $N$ is a
    complex of discrete $\Lambda$-modules.
    
    Now $\Lambda/\mathfrak M^l$ is such a finite module, hence
    \[\Rder\Gamma_{\underline{\mathfrak M}}(\Lambda) =
      [-d-s]\circ \varinjlim_l (\Rder\Gamma_{\underline{\mathfrak
          m}}(\Lambda/\mathfrak M^l)(\chi))^\vee =
      \Lambda(\chi)^\vee[-r].\]The proposition now follows at once if
    we observe that $\Lambda\cong\Lambda(\chi)$ as a $\Lambda$-module
    via $g\mapsto \chi(g)g$.
  \end{proof}
\end{proposition}

\begin{theorem}[Local duality for Iwasawa
  algebras]\label{thm:local-duality-for-iw-alg}
  \[\Rder\Gamma_{\underline{\mathfrak M}} \cong [-r] \circ \Pi \circ
  \RHom_{\Lambda}(-,\Lambda)\] on $\Dder^b_c(\Lambda)$.
  \begin{proof}
    Because of
    \cref{prop:g-poincare-then-top-loc-coh-is-pont-dual-for-ring}, this
    follows verbatim as in \cite[theorem 6.3]{MR0224620}: The functors
    $\Rder^r\Gamma_{\underline{\mathfrak M}}$ and
    $\Hom_\Lambda(-,\Lambda)^\vee$ are related by a pairing of
    $\Ext$-groups, are both covariant and right-exact and agree on
    $\Lambda$, hence also agree on finitely generated modules. As the
    complex $\Rder\Gamma_{\underline{\mathfrak M}}(\Lambda)$ is
    concentrated in degree $r$, the same argument as in
    \cref{prop:dual-of-iw-adj-is-twist-of-local-coh-of-coeff} shows
    that the left derivation of
    $\Rder^r\Gamma_{\underline{\mathfrak M}}$ is just
    $[r]\circ \Rder\Gamma_{\underline{\mathfrak M}}$ and the result
    follows.
  \end{proof}
\end{theorem}

\section{Torsion in Iwasawa Cohomology}

There are notions of both local and global Iwasawa cohomology. Our result about
their torsion below holds in both cases and we will first deal with the local
case.

In both subsections, $R$ is a commutative Noetherian pro-$p$ local ring of
residue characteristic $p$.

\subsection{Torsion in Local Iwasawa Cohomology}

Let $K$ be a finite extension of $\Q_p$ and $K_\infty|K$ a Galois extension with
an analytic pro-$p$ Galois group $G$ without elements of finite order. Let $T$
be a finitely generated $\Lambda=R[[G]]$-module and set set $A=T\otimes_R
R^\vee$. Due to Lim and Sharifi we have the following spectral sequence
(stemming from an isomorphism of complexes in the derived category). Write
$H^i_{\mathrm{Iw}}(K_\infty, T) = \varprojlim_{K'} H^{i}(G_{K'}, T)$ where the
limit is taken with respect to the corestriction maps over all finite field
extensions $K'|K$ contained in $K_\infty.$

\begin{theorem}\label{thm:lim-sharifi-spec-seq-local}
  There is a convergent spectral sequence \[\E^i_\Lambda(H^j(G_K, A)^\vee)
    \Longrightarrow H^{i+j}_{\mathrm{Iw}}(K,T).\]
  \begin{proof}
    This is \cite[theorem 4.2.2, remark 4.2.3]{MR3084561}.
  \end{proof}
\end{theorem}

\begin{theorem}\label{thm:torsion-in-local-iw-coh-for-poinc-grps}
  If $G$ is a pro-$p$ Poincaré group of dimension $s\geq 2$ with dualizing
  character $\chi\colon G\ra\Z_p^\times$ and if $R$ is regular,
  then \[\tor_\Lambda H^1_{\mathrm{Iw}}(K_\infty, T) = 0.\] If $s=1$, then
  \[\tor_\Lambda H^1_{\mathrm{Iw}}(K_\infty, T) \cong \Hom_R(
    (T^*)_G, R)(\chi^{-1}),\] where $T^* = \Hom_R(T,R)$.
  \begin{proof}
    The exact five-term sequence attached to
    \cref{thm:lim-sharifi-spec-seq-local} starts like this:
    \[ \begin{tikzcd} 0 \arrow{r} & \E^1_\Lambda(H^0(G_K,A)^\vee) \arrow{r} &
        H^1_{\mathrm{Iw}}(K_\infty, T) \arrow{r} & \E^0_\Lambda(H^1(G_K,A)^\vee)
        \arrow{d}\\ & & & \E^2_\Lambda(H^0(G_K,A)^\vee) \end{tikzcd}\] Note that
    $\E^2_\Lambda(M)$ is pseudo-null and hence $\Lambda$-torsion for every
    finitely generated module $M$, which follows from the spectral sequence
    attached to the isomorphism $\RHom_\Lambda(-,\Lambda)\circ
    \RHom_\Lambda(-,\Lambda) \cong \id$ on $\Dder^b_c(\RMod{\Lambda})$.
    Furthermore $\E^0_\Lambda(M)$ is always $\Lambda$-torsion free, as $\Lambda$
    is integral. It follows that $\tor_\Lambda H^1_{\mathrm{Iw}}(K_\infty, T)
    \subseteq \E^1_\Lambda(H^0(G_K,A)^\vee)$. As the latter is
    $\Lambda$-torsion,
    \[\tor_\Lambda H^1_{\mathrm{Iw}}(K_\infty, T) =
      \E^1_\Lambda(H^0(G_K,A)^\vee).\] The result now follows
    immediately from
    \cref{prop:explicit-iw-adj}
  \end{proof}
  
\end{theorem}

\subsection{Torsion in Global Iwasawa Cohomology}

Let $K$ be a finite extension of $\Q$ and $S$ a finite set of places of $K$. Let
$K_S$ be the maximal extension of $K$ which is unramified outside $S$ and
$K_\infty|K$ a Galois extension contained in $K_S$. Suppose that
$G=G(K_\infty|K)$ is an analytic pro-$p$ group without elements of finite order.
Let $T$ be a finitely generated $\Lambda=R[[G]]$-module and set $A=T\otimes_R
R^\vee$. Due to Lim and Sharifi we have the following spectral sequence
(stemming from an isomorphism of complexes in the derived category). Write
$H^i_{\mathrm{Iw}}(K_\infty, T) = \varprojlim_{K'} H^{i}(G(K_S|K'), T)$ where
the limit is taken with respect to the corestriction maps over all finite field
extensions $K'|K$ contained in $K_\infty.$

\begin{theorem}\label{thm:lim-sharifi-spec-seq-global}
  There is a convergent spectral sequence
  \[\E^i_\Lambda(H^j(G_S, A)^\vee) \Longrightarrow H^{i+j}_{\mathrm{Iw}}(K_\infty,
    T).\]
\begin{proof}
    This follows from \cite[theorem 4.5.1]{MR3084561}.
\end{proof}
\end{theorem}

From this we derive the following.

\begin{theorem}\label{thm:torsion-in-global-iw-coh-for-poinc-grps}
  If $G$ is a pro-$p$ Poincaré group of dimension $s\geq 2$ with
  dualizing character $\chi\colon G\ra\Z_p^\times$ and if $R$ is
  regular, then \[\tor_\Lambda H^1_{\mathrm{Iw}}(K_\infty, T) = 0.\]
  If $s=1$, then
  \[\tor_\Lambda H^1_{\mathrm{Iw}}(K_\infty, T) \cong \Hom_R(
    (T^*)_G, R)(\chi^{-1}),\]
  where $T^* = \Hom_R(T,R)$.
  \begin{proof}
    Replace “$G_K$” with “$G_S$” in the proof of \cref{thm:torsion-in-local-iw-coh-for-poinc-grps}.
  \end{proof}
\end{theorem}

\bibliographystyle{hamsalpha.bst}
\bibliography{bib.bib}{}

\end{document}